\documentclass[12pt]{amsart}
\usepackage{amssymb}

\usepackage{amsfonts}

\newtheorem{theorem}{Theorem}[section]
\theoremstyle{plain}

\newtheorem{proposition}[theorem]{Proposition}

\numberwithin{equation}{section}
\input{tcilatex}

\begin{document}
\title[Comparison Theorem]{Comparison Theorem and Geometric Realization of
Representations. }
\author{Tim Bratten}
\address{Facultad de Ciencias Exactas\\
UNICEN\\
Paraje Arroyo Seco\\
Tandil, Argentina}
\email{bratten@exa.unicen.edu.ar}

\begin{abstract}
In this paper we generalize the comparison theorem of Hecht and Taylor to
arbitrary parabolic subalgebras of a complex reductive Lie algebra and then
apply our generalized comparison theorem to obtain results about the
geometric realization of representations in flag spaces.
\end{abstract}

\maketitle

\section{Introduction}

This manuscript concerns a basic homological property of representations for
a reductive Lie group, called the comparison theorem, and the relation this
property has to the realization of representations in complex flag spaces.

Historically speaking, the realization of representations in complex flag
spaces has played a central role in the theory. One of the earliest
constructions, the parabolic induction, defines representations as the
sections of homogeneous vector bundles defined over certain closed orbits in
complex flag spaces. Schmid's realization of the discrete series \cite%
{thesis}, generalizing the Borel-Weil-Bott theorem, gave a defining
alternative to the parabolic induction, finding the missing representations
on the sheaf cohomology groups of certain homogeneous holomorphic line
bundles defined on open orbits in full flag spaces.

The problem of understanding sheaf cohomologies of homogeneous holomorphic
vector bundles turned out to be somewhat delicate, and it took some time
until general results were obtained, but meanwhile the localization theory
of Beilinson and Bernstein \cite{algebraic} provided a sort of universal
geometric realization, defined in the full flag space, for irreducible
Harish-Chandra modules. Via localization, many irreducible Harish-Chandra
modules are nicely realized as certain standard geometric objects (a precise
criteria for this is known \cite{irreducible}) but in general it seems quite
difficult to understand the localization of an irreducible representation.
The analytic localization theory of Hecht and Taylor \cite{analytic} gives a
global counterpart to the Beilinson-Bernstein algebraic theory. A main
result of the analytic theory shows that the compactly supported cohomology
of the polarized sections of an irreducible homogeneous vector bundle
realizes the minimal globalization of the cohomology of an associated
standard Beilinson-Bernstein sheaf.

Although the Hecht-Taylor result realizes, for example, all of the tempered
representations, many irreducible representations are not realized as
standard modules in a full flag space. Thus one considers analogous
constructions defined on arbitrary flag spaces. Along these lines, Wong \cite%
{wong} studied the representations obtained on the sheaf cohomologies of
finite rank homogeneous holomorphic vector bundles defined over certain open
orbits in generalized flag spaces, proving a special case of a conjecture by
Vogan \cite{vogan}. Via the methods of algebraic and analytic localization
in flag spaces, the author of this study was able to prove the general
version of Vogan's conjecture and give a realization for all of the
Harish-Chandra modules defined by cohomological parabolic induction \cite%
{realizing}.

In this study we consider how the mesh of localization theories works, in
complete generality, for complex flag spaces. In particular, we establish
the full extent to which previous results about the realization of standard
modules can generalize. As an intermediate result we obtain Theorem 4.1,
which generalizes the Hecht-Taylor comparison theorem \cite{comparison1} to
arbitrary orbits in flag spaces, provided one makes a finiteness assumption.
A main result, Theorem 5.5, applies the generalized comparison theorem to
show that the Hecht-Taylor realization of standard modules extends naturally
to a class of orbits we call affinely oriented (this includes all open
orbits, and therefore all the homogeneous holomorphic vector bundles). We
finish our study by analyzing an example showing how the situation works in
the case where the orbit in question is not affinely oriented. In
particular, we show that Theorem 5.5 fails to hold.

This study is organized as follows. In Section 2 we introduce the algebraic
and analytic localization in flag spaces and establish some simple results
we will use. In Section 3 we define the standard modules. In Section 4 we
prove the comparison theorem and in Section 5 we establish our main result.
In the last section we consider the $SU(n,1)$ action in complex projective
space and see how the main result fails when the orbit is not affinely
oriented. 

\section{Algebraic and Analytic Localization}

In this section we introduce the minimal globalization, define the flag
spaces, consider the generalized TDOs and establish some facts about the
algebraic and analytic localization in complex flag spaces.

Throughout this study, $G_{0}$ denotes a reductive group of Harish-Chandra
class with Lie algebra $\mathfrak{g}_{0}$ and $\mathfrak{g}$ denotes the
complexification of $\mathfrak{g}_{0}$. We fix a maximal compact subgroup $%
K_{0}$ of $G_{0}$ and let $K$ denote the complexification of $K_{0}$. $G$
indicates the complex adjoint group of $\mathfrak{g}$.

\bigskip

\noindent \textbf{Minimal Globalization.} \ By definition, a \emph{\
Harish-Chandra module }is a finite-length $\mathfrak{g}-$module equipped
with a compatible, algebraic $K-$action. For example, the set of $K_{0}-$%
finite vectors in an irreducible unitary representation for $G_{0}$ is a
Harish-Chandra module.

Let $M$ be a Harish-Chandra module. A \emph{globalization} of $M$ means a
finite-length, admissible representation for $G_{0}$ in a complete, locally
convex space whose underlying space of $K_{0}-$finite vectors is $M$. By now
there are known to exist several canonical and functorial globalizations of
Harish-Chandra modules, including the remarkable \emph{minimal globalization}%
, whose existence was first proved by Schmid \cite{minimal}. The minimal
globalization is functorial and embeds continuously and $G_{0}-$%
equivariantly in any corresponding globalization. Indeed, as a functor the
minimal globalization is exact and surjects onto the space of analytic
vectors in a Banach space globalization \cite{minimal2}.

\bigskip

\noindent \textbf{Flag Spaces.} \ By a\emph{\ complex flag space} for $G_{0}$
we mean a complex projective homogeneous $G-$space $Y$. The complex flag
spaces are constructed as follows. By definition, a \emph{Borel subalgebra}
of $\mathfrak{g}$ is a maximal solvable subalgebra. A basic fact is that $G$
acts transitively on the set of Borel subalgebras and the resulting
homogeneous space is a complex projective variety $X$ called\emph{\ the full
flag space} for $G_{0}$. A complex subalgebra that contains a Borel
subalgebra is called a \emph{parabolic subalgebra} of $\mathfrak{g}$. If one
fixes a Borel subalgebra $\mathfrak{b}$ of $\mathfrak{g}$, then each
parabolic subalgebra is $G-$conjugate to a unique parabolic subalgebra
containing $\mathfrak{b}$. The resulting space $Y$ of $G-$conjugates to a
given parabolic subalgebra is a complex flag space and each complex flag
space is realized this way.

Let $X$ denote the full flag space and suppose $Y$ is a complex flag space.
For $x\in X$ and $y\in Y$ we let $\mathfrak{b}_{x}$ and $\mathfrak{p}_{y}$
indicate, respectively, the corresponding Borel and parabolic subalgebras of 
$\mathfrak{g}$. From the above discussion it follows that there exits a
unique $G-$equivariant projection 
\begin{equation*}
\pi :X\rightarrow Y
\end{equation*}%
given by $\pi (x)=y$ where $y\in Y$ is the unique point such that $\mathfrak{%
b}_{x}\subseteq \mathfrak{p}_{y}$. $\pi $ is called \emph{the natural
projection}.

We will need to treat $Y$ as both a complex analytic manifold with its sheaf
of holomorphic functions $\mathcal{O}_{Y}$ and as an algebraic variety $Y^{%
\text{alg}}$ with the Zariski topology and the corresponding sheaf of
regular functions $\mathcal{O}_{Y^{\text{alg}}}$.

\bigskip

\noindent \textbf{Generalized Sheaves of TDOs}. \textbf{\ }Let $U(\mathfrak{g%
})$ denote the enveloping algebra of $\mathfrak{g}$ and let $Z(\mathfrak{g})$
be the center of $U(\mathfrak{g})$. By definition, a $\mathfrak{g}-$\emph{%
infinitesimal character} is a homomorphism of algebras 
\begin{equation*}
\Theta :Z(\mathfrak{g})\rightarrow \mathbb{C}\text{.}
\end{equation*}%
By a fundamental result of Harish-Chandra, one can parametrize the $%
\mathfrak{g}-$infinitesimal characters as follows. Let $\mathfrak{h}^{\ast }$
be the \emph{Cartan dual }for $\mathfrak{g}$ (definitions as in \cite%
{realizing}, Section 2 and Section 3). There is a naturally defined set of 
\emph{roots }$\Sigma \subseteq \mathfrak{h}^{\ast }$ for $\mathfrak{h}$ in $%
\mathfrak{g}$ and a corresponding subset of \emph{positive roots} $\Sigma
^{+}\subseteq \Sigma $. Let $W$ denote the Weyl group for $\mathfrak{h}%
^{\ast }$ induced by the roots of $\mathfrak{h}$ in $\mathfrak{g}$. Then
there is a natural 1-1 correspondence between the set of $\mathfrak{g}-$%
infinitesimal characters and the quotient 
\begin{equation*}
\mathfrak{h}^{\ast }/W\text{.}
\end{equation*}%
Given $\lambda \in \mathfrak{h}^{\ast }$ and an infinitesimal character $%
\Theta $, we write $\lambda \in \Theta $ to indicate that $\Theta $
corresponds to the Weyl group orbit $W\cdot \lambda $ in the given
parametrization. $\Theta $ is called\emph{\ regular }when the corresponding
Weyl group orbit has the order of $W$ elements. $\lambda \in \mathfrak{h}%
^{\ast }$ is called \emph{antidominant} if 
\begin{equation*}
\overset{\vee }{\alpha }(\lambda )\notin \left\{ 1,2,3,\ldots \right\} \text{
for each positive root }\alpha \in \Sigma ^{+}\text{ .}
\end{equation*}%
We also introduce the following notation: given a $\mathfrak{g}-$%
infinitesimal character $\Theta $ we let $U_{\Theta }$ denote the algebra
obtained as the quotient of the enveloping algebra $U(\mathfrak{g})$ by the
ideal generated from the ideal in $Z(\mathfrak{g})$ corresponding to $\Theta 
$. In particular, a $U_{\Theta }$-module is just a $\mathfrak{g}-$module
with infinitesimal character $\Theta $.

To each $\lambda \in \mathfrak{h}^{\ast }$, Beilinson and Bernstein
associate a twisted sheaf of differential operators (TDO) $\mathcal{D}%
_{\lambda }^{\text{alg}}$ defined on the algebraic variety $X^{\text{alg}}$ 
\cite{algebraic}. In our parametrization, the sheaf of differential
operators on $X^{\text{alg}}$ is $\mathcal{D}_{-\rho }^{\text{alg}}$, where $%
\rho \in \Sigma ^{+}$ is one half the sum of the positive roots. Beilinson
and Bernstein prove that 
\begin{equation*}
H^{p}(X,\mathcal{D}_{\lambda }^{\text{alg}})=0\;\text{for }p>0\text{ and \
that }U_{\Theta }\cong \Gamma (X,\mathcal{D}_{\lambda }^{\text{alg}})\text{
\ }
\end{equation*}%
where $\Theta =W\cdot \lambda $. In particular 
\begin{equation*}
\Gamma (X,\mathcal{D}_{\lambda }^{\text{alg}})\cong \Gamma (X,\mathcal{D}%
_{w\lambda }^{\text{alg}})\text{ \ for }w\in W\text{.}
\end{equation*}

Let $\pi _{\ast }$ denote the direct image in the category of sheaves. We
consider the \emph{generalized sheaf of TDOs} $\pi _{\ast }(\mathcal{D}%
_{\lambda }^{\text{alg}})$ defined on $Y^{\text{alg}}$. Observe that 
\begin{equation*}
U_{\Theta }=\Gamma (Y,\pi _{\ast }(\mathcal{D}_{\lambda }^{\text{alg}}))
\end{equation*}%
where $\Theta =W\cdot \lambda $.

Given $y\in Y$ let $\mathfrak{p}_{y}$ be the corresponding parabolic
subalgebra of $\mathfrak{g}$ and let $\mathfrak{u}_{y}$ denote nilradical of 
$\mathfrak{p}_{y}$. The \emph{Levi quotient is }given by \emph{\ } \emph{\ } 
\begin{equation*}
\mathfrak{l}_{y}=\mathfrak{p}_{y}/\mathfrak{u}_{y}\text{.}
\end{equation*}%
Since Cartan subalgebras of $\mathfrak{g}$ contained in $\mathfrak{p}_{y}$
are naturally identified with Cartan subalgebras of $\mathfrak{l}_{y}$ one
can, in natural way, identify $\mathfrak{h}^{\ast }$ with the Cartan dual
for the reductive Lie algebra $\mathfrak{l}_{y}$. This induces a set of
roots 
\begin{equation*}
\Sigma _{Y}\subseteq \Sigma
\end{equation*}%
of $\mathfrak{h}^{\ast }$ in $\mathfrak{l}_{y}$, and a Weyl group $%
W_{Y}\subseteq W$, generated by reflections coming from the elements of $%
\Sigma _{Y}$. As suggested in the notation, these subsets are independent of
the point $y$. One proves that 
\begin{equation*}
\pi _{\ast }(\mathcal{D}_{\lambda }^{\text{alg}})\cong \pi _{\ast }(\mathcal{%
D}_{w\lambda }^{\text{alg}})\text{ \ for }w\in W_{Y}\text{. }
\end{equation*}%
We say that $\lambda \in \mathfrak{h}^{\ast }$ is \emph{antidominant for }$Y$
if $\lambda $ is $W_{Y}-$conjugate to an antidominant element of $\mathfrak{h%
}^{\ast }$. Generalizing the result of Beilinson and Bernstein for the
twisted sheaves of differential operators on $X^{\text{alg}}$, it has been
show \cite{chang} that if $\mathcal{F}$ is a quasicoherent sheaf of $\pi
_{\ast }(\mathcal{D}_{\lambda }^{\text{alg}})-$modules on $Y^{\text{alg}}$
and if $\lambda $ is antidominant for $Y$ then 
\begin{equation*}
H^{p}(Y,\mathcal{F})=0\text{ \ for }p>0\text{.}
\end{equation*}

\bigskip

\noindent \textbf{Algebraic and Analytic Localization.} \ Given a $\mathfrak{%
g}-$module $M$ with infinitesimal character $\Theta $ and a choice of $%
\lambda \in \Theta $ we define the (algebraic) localization of $M$ to $Y^{%
\text{alg}}$ as the sheaf of $\pi _{\ast }(\mathcal{D}_{\lambda }^{\text{alg}%
})-$modules given by 
\begin{equation*}
\Delta _{\lambda }^{\text{alg}}(M)=M\otimes _{U_{\Theta }}\pi _{\ast }(%
\mathcal{D}_{\lambda }^{\text{alg}})\text{.}
\end{equation*}%
Thus 
\begin{equation*}
\Delta _{\lambda }^{\text{alg}}(M)\cong \Delta _{w\lambda }^{\text{alg}}(M)%
\text{ \ for }w\in W_{Y}\text{.}
\end{equation*}%
Generalizing the Beilinson-Bernstein result for a full flag space, it has
been shown that when $\Theta $ is a regular infinitesimal character and $%
\lambda \in \Theta $ is antidominant for $Y$ then the localization functor
and the global sections on $Y^{\text{alg}}$ define an equivalence of
categories between the category of $U_{\Theta }-$modules and the category of
quasicoherent $\pi _{\ast }(\mathcal{D}_{\lambda }^{\text{alg}})-$modules 
\cite{chang}.

While the algebraic localization functor yields interesting results when
applied to Harish-Chandra modules, the analytic localization of Hecht and
Taylor allows one to study the geometric realization for the minimal
globalization. Since the analytic localization takes into account the
topology of a module, we introduce a few relevant concepts. By definition, a%
\emph{\ DNF space} is topological vector space whose strong dual is a
nuclear Fr\'{e}chet space. A \emph{DNF} $U_{\Theta }$\emph{-module }is a $%
\mathfrak{g}-$module with infinitesimal character $\Theta $, defined on a
DNF-space $M$ such that the corresponding linear operators 
\begin{equation*}
m\mapsto \xi \cdot m\text{ \ for }m\in M\text{ and }\xi \in \mathfrak{g}
\end{equation*}%
are continuous. Observe that a finitely generated $U_{\Theta }-$module has a
unique DNF topology, when considered as a topological direct sum of finite
dimensional subspaces.

Given $\lambda \in \mathfrak{h}^{\ast }$, let $\mathcal{D}_{\lambda }$
denote the corresponding TDO with holomorphic coefficients defined on the
complex manifold $X$. We consider the generalized TDO $\pi _{\ast }(\mathcal{%
D}_{\lambda })$ defined on the complex flag space $Y$. Let $\mathcal{O}_{Y}$
denote the sheaf of holomorphic functions on $Y$. Since the sheaf $\pi
_{\ast }(\mathcal{D}_{\lambda })$ is locally free as a sheaf of $\mathcal{O}%
_{Y}-$modules with countable geometric fiber, there is a natural DNF
topology defined on the space of sections of $\pi _{\ast }(\mathcal{D}%
_{\lambda })$ over compact subsets of $Y$ \cite{analytic}. When $M$ is a DNF 
$U_{\lambda }$-module then, using the completed tensor product, one can
define a sheaf 
\begin{equation*}
\Delta _{\lambda }(M)=\pi _{\ast }(\mathcal{D}_{\lambda })\overset{\wedge }{%
\otimes }_{U_{\lambda }}M
\end{equation*}%
of $\pi _{\ast }(\mathcal{D}_{\lambda })-$modules carrying a natural
topological vector space structure defined over compact subsets of $Y$. In
case the induced topologies on the geometric fibers of $\Delta _{\lambda
}(M) $ are Hausdorff, then $\Delta _{\lambda }(M)$ is a DNF sheaf of $\pi
_{\ast }(\mathcal{D}_{\lambda })-$modules \cite{realizing}.

Given a $U_{\Theta }-$module $M$, the Hochschild resolution $F_{\cdot }(M)$
of $M$ is the canonical resolution of $M$ by free $U_{\Theta }-$modules
where 
\begin{equation*}
F_{p}(M)=\otimes ^{p+1}U_{\Theta }\otimes M.
\end{equation*}
When $M$ is a DNF $U_{\Theta }-$module, then $F_{\cdot }(M)$ is a complex of
DNF $U_{\Theta }$-modules and $\Delta _{\lambda }(F_{\cdot }(M))$ is a
functorially defined complex of DNF sheaves of $\pi _{\ast }(\mathcal{D}%
_{\lambda })-$modules called \emph{the analytic localization of} $M$ \emph{to%
} $Y$ \emph{with respect to} \emph{the parameter} $\lambda \in \Theta
=W\cdot \lambda $. Observe that 
\begin{equation*}
\Delta _{\lambda }(F_{\cdot }(M))\cong \Delta _{w\lambda }(F_{\cdot }(M))%
\text{ \ for }w\in W_{Y}\text{.}
\end{equation*}%
It is not hard to show that this complex of sheaves has hypercohomology
naturally isomorphic to $M$. We shall analyze in more detail the results of
the analytic localization as applied to minimal globalizations as our study
advances, but right now we want to point out that when $M$ is a minimal
globalization then $G_{0}$ acts naturally on the homology sheaves of $\Delta
_{\lambda }(F_{\cdot }(M))$ . In particular, $G_{0}$ acts on $F_{\cdot }(M)$
by the tensor product of the adjoint action with the action on $M$. Although
this action is not compatible with the left $\mathfrak{g}-$action, the two
actions are homotopic. Coupling the $G_{0}-$action on $\pi _{\ast }(\mathcal{%
D}_{\lambda })$ with the $G_{0}-$action on $F_{\cdot }(M)$, one obtains a $%
G_{0}-$action on $\Delta _{\lambda }(F_{\cdot }(M))$. \ 

\bigskip

\noindent \textbf{Localization and Geometric Fibers. \ }In order to prove
the comparison theorem in Section 4, we will use some simple facts about the
localization functors and geometric fibers, which we summarize in the
following two propositions. The first proposition says that, for computing
geometric fibers, the algebraic and analytic localizations yield the same
result. In particular, consider the sheaves $\mathcal{O}_{Y}$ on $Y$ and $%
\mathcal{O}_{Y^{\text{alg}}}$ on $Y^{\text{alg}}$. Fix $y\in Y$. If $%
\mathcal{F}$ is sheaf of $\mathcal{O}_{Y}-$modules on $Y$ and $\mathcal{H}$
is a sheaf of $\mathcal{O}_{Y^{\text{alg}}}-$modules on $Y^{\text{alg}}$,
put 
\begin{equation*}
T_{y}(\mathcal{F})=\mathbb{C\otimes }_{\mathcal{O}_{Y}\mid _{y}}\mathcal{F}%
\text{ \ \ and \ }T_{y}^{\text{alg}}(\mathcal{F})=\mathbb{C\otimes }_{%
\mathcal{O}_{Y^{\text{alg}}}\mid _{y}}\mathcal{H}
\end{equation*}%
where $\mathcal{O}_{Y}\mid _{y}$ and $\mathcal{O}_{Y^{\text{alg}}}\mid _{y}$
denote the respective stalks of $\mathcal{O}_{Y}$ and $\mathcal{O}_{Y^{\text{%
alg}}}$ over the point $y$. When $\mathcal{F}$ is a sheaf $\mathcal{D}%
_{\lambda }$-modules (respectively a sheaf of $\mathcal{D}_{\lambda }^{\text{%
alg}}$-modules) then, in a natural way, $T_{y}(\mathcal{F})$ (respectively $%
T_{y}^{\text{alg}}(\mathcal{F})$), is a module for the corresponding Levi
quotient.

\begin{proposition}
Let $M$ be a DNF $U_{\Theta }$-module and choose $\lambda \in \Theta $. Let $%
Y$ be a complex flag space and choose $y\in Y$. Let $\mathfrak{l}_{y}$
denote the corresponding Levi quotient. Then there is a natural isomorphism 
\begin{equation*}
T_{y}\circ \Delta _{\lambda }(F_{\bullet }(M))\cong T_{y}^{\text{alg}}\circ
\Delta _{\lambda }^{\text{alg}}(F_{\bullet }(M))\text{ }
\end{equation*}
of complexes of $\mathfrak{l}_{y}$-modules.
\end{proposition}

\medskip

\noindent \textbf{Proof:} \ Via the natural inclusion 
\begin{equation*}
\pi _{\ast }(\mathcal{D}_{\lambda }^{\text{alg}})\mid _{y}\rightarrow \pi
_{\ast }(\mathcal{D}_{\lambda })\mid _{y}
\end{equation*}%
one obtains an isomorphism 
\begin{equation*}
T_{y}^{\text{alg}}(\pi _{\ast }(\mathcal{D}_{\lambda }^{\text{alg}}))\cong
T_{y}(\pi _{\ast }(\mathcal{D}_{\lambda }))
\end{equation*}%
of left $\mathfrak{l}_{y}-$modules. Thus there is a corresponding a natural
isomorphism 
\begin{equation*}
T_{y}^{\text{alg}}(\Delta _{\lambda }^{\text{alg}}(N))\cong T_{y}(\pi _{\ast
}(\mathcal{D}_{\lambda }))\otimes _{U_{\Theta }}N
\end{equation*}%
where $N$ is a $U_{\Theta }-$module. On the hand, since $T_{y}(\pi _{\ast }(%
\mathcal{D}_{\lambda }))$ has countable dimension, it follows that the
natural inclusion determines an isomorphism 
\begin{equation*}
T_{y}(\pi _{\ast }(\mathcal{D}_{\lambda }))\otimes _{U_{\Theta }}N\cong
T_{y}\circ \Delta _{\lambda }(N)
\end{equation*}%
when $N$ is DNF $U_{\Theta }-$module. $\blacksquare $

\bigskip

We will also use the following base change formulas. Let 
\begin{equation*}
X_{y}=\pi ^{-1}(\left\{ y\right\}
\end{equation*}%
be the fiber in $X$ over $y$ and let 
\begin{equation*}
i:X_{y}\rightarrow X
\end{equation*}%
denote the inclusion. Suppose $i^{-1}$ denotes the corresponding inverse
image (in this case: the restriction) in the category of sheaves. If $%
\mathcal{F}$ is sheaf of $\mathcal{O}_{X}-$modules on $X$ and $\mathcal{H}$
is a sheaf of $\mathcal{O}_{X^{\text{alg}}}$-modules on $X^{\text{alg}}$, we
put 
\begin{equation*}
i^{\ast }(\mathcal{F})=\mathcal{O}_{X_{y}}\otimes _{i^{-1}(\mathcal{O}%
_{X})}i^{-l}(\mathcal{F})\text{ }
\end{equation*}%
and 
\begin{equation*}
i_{\text{alg}}^{\ast }(\mathcal{H})=\mathcal{O}_{X_{y}^{\text{alg}}}\otimes
_{i^{-1}(\mathcal{O}_{X^{\text{alg}}})}i^{-l}(\mathcal{H}).
\end{equation*}%
When $\mathcal{F}$ is a sheaf $\mathcal{D}_{\lambda }-$modules (respectively
a sheaf of $\mathcal{D}_{\lambda }^{\text{alg}}-$modules) then, in a natural
way, $\Gamma (X_{y},i^{\ast }(\mathcal{F}))$ (respectively $\Gamma (X_{y},i_{%
\text{alg}}^{\ast }(\mathcal{F}))$), is a module for the corresponding Levi
quotient. The base change formulas are the following two results.

\begin{proposition}
Let $M$ be a DNF $U_{\Theta }$-module and let $N$ be a $U_{\Theta }$-module.
Choose $\lambda \in \Theta $. Let $Y$ be a complex flag space and choose $%
y\in Y$. Let $\Delta _{X,\lambda }$ and $\Delta _{Y,\lambda }$ denote the
corresponding analytic localizations to $X$ and $Y$ and let $\Delta
_{X,\lambda }^{\text{alg}}$ and $\Delta _{Y,\lambda }^{\text{alg}}$ denote
the corresponding algebraic localizations to $X^{\text{alg}}$ and $Y^{\text{%
alg}}$. Then, using the above notations, we have the following natural
isomorphisms of complexes of $\mathfrak{l}_{y}$-modules: 
\begin{equation*}
\text{(a) }T_{y}\circ \Delta _{Y,\lambda }(F_{\bullet }(M))\cong \Gamma
(X_{y},i^{\star }\circ \Delta _{X,\lambda }(F_{\bullet }(M))\text{;}
\end{equation*}%
\begin{equation*}
\text{(b) }T_{y}^{\text{alg}}\circ \Delta _{Y,\lambda }^{\text{alg}%
}(F_{\bullet }(N))\text{ }\cong \Gamma (X_{y},i_{\text{alg}}^{\ast }\circ
\Delta _{X,\lambda }^{\text{alg}}(F_{\bullet }(N)).\text{ }
\end{equation*}
\end{proposition}

\noindent \textbf{Proof:} Equation (a) is shown in \cite{realizing},
Proposition 3.3 and Equation (b) can be proved in exactly the same way. $%
\blacksquare $

\section{Standard Modules in Flag Spaces}

In this section we review the Matsuki duality, consider polarized
representations\ for the stabilizer and introduce the standard modules. We
finish the section by summarizing a result, due to Hecht and Taylor, that
characterizes the analytic localization of a minimal globalization to the
full flag space. \ 

\bigskip

\noindent \textbf{Matsuki Duality.} \ \ Let $Y$ be a complex flag space. It
is known that $G_{0}$ acts with finitely many orbits in a $Y$. We will need
to use the following geometric property relating the $G_{0}$ and $K-$actions
in $Y$, referred to as \emph{Matsuki duality}. Let 
\begin{equation*}
\theta :\mathfrak{g}\rightarrow \mathfrak{g}
\end{equation*}%
denote the complexified Cartan involution arising from $K_{0}$ and let 
\begin{equation*}
\tau :\mathfrak{g}\rightarrow \mathfrak{g}
\end{equation*}%
denote the conjugation corresponding to $\mathfrak{g}_{0}$. A subalgebra of $%
\mathfrak{g}$ is called \emph{stable} if it is invariant under both $\theta $
and $\tau $. A point $y\in Y$ is called \emph{special} if $\mathfrak{p}_{y}$
contains a stable Cartan subalgebra of $\mathfrak{g}$. A $G_{0}-$orbit $S$
is said to be \emph{Matsuki dual} to a $K-$orbit $Q$ when $S\cap Q$ contains
a special point. Since it is known that the set of special points in a $%
G_{0}-$orbit, or in a $K-$orbit, forms a nonempty $K_{0}-$homogeneous
submanifold it follows that Matsuki duality gives a 1-1 correspondence
between the $G_{0}-$orbits and the $K-$orbits on $Y$ \cite{matsuki}.

\bigskip

\noindent \textbf{Polarized Modules}. \ Suppose $y\in Y$ and let $G_{0}[y]$
denote the stabilizer of $y$ in $G_{0}$. Let 
\begin{equation*}
\omega :G_{0}[y]\rightarrow GL(V)
\end{equation*}%
be a representation in a finite-dimensional complex vector space $V$. A
compatible, linear $\mathfrak{p}_{y}-$action in $V$ is called a \emph{%
polarization} if the nilradical $\mathfrak{u}_{y}$ acts trivially. In other
words: a polarized $G_{0}[y]-$module is a nothing but a finite-dimensional $(%
\mathfrak{l}_{y},G_{0}[y])-$module. In case $\mathfrak{p}_{y}$ contains a%
\emph{\ real Levi factor} (that is: a complementary subalgebra to the
nilradical that is invariant under $\tau $) then an irreducible
representation always has a unique polarization, but in general compatible $%
\mathfrak{p}_{y}-$actions need not exist. For example, suppose $\mathfrak{c}$
is a stable Cartan subalgebra of $\mathfrak{g}$, $\mathfrak{b}$ is a Borel
subalgebra containing $\mathfrak{c}$ and $\alpha \in \mathfrak{c}^{\ast }$
is a simple root of $\mathfrak{c}$ in $\mathfrak{b}$. Let $\mathfrak{g}%
^{\alpha }$ and $\mathfrak{g}^{-\alpha }$ denote the corresponding root
subspaces of $\mathfrak{c}$ in $\mathfrak{g}$ and define 
\begin{equation*}
\mathfrak{p}_{y}=\mathfrak{g}^{-\alpha }+\mathfrak{b}\text{.}
\end{equation*}%
Then $\mathfrak{p}$ is a parabolic subalgebra of $\mathfrak{g}$. Assume that
the root $\alpha $ is complex, that is: 
\begin{equation*}
\mathfrak{p}\cap \tau (\mathfrak{p})=\mathfrak{c}
\end{equation*}%
and let $C_{0}$ be the Cartan subgroup of $G_{0}$ corresponding to $%
\mathfrak{c}.$ Then the character of $C_{0}$ given by the adjoint action of $%
C_{0}$ in $\mathfrak{g}^{\alpha }$ extends uniquely to a character of $%
G_{0}[y]$ [cite], but there is no associated polarization.

Even though polarizations need not exist, they are unique when they do
exist. In particular, suppose $V$ is a $G_{0}[y]$-module with two
polarizations. Then there are two $\mathfrak{l}_{y}$-actions in $V$ that
coincide on a parabolic subalgebra of $\mathfrak{l}_{y}$ [cite]. From the
theory of finite-dimensional $\mathfrak{l}_{y}$-modules, it follows that the
two $\mathfrak{l}_{y}$-actions are identical.

On the other hand, if $V$ is a finite-dimensional irreducible $(\mathfrak{p}%
_{y},G_{0}[y])-$module then the $\mathfrak{p}_{y}-$action is necessarily a
polarization, since the subspace of vectors annihilated by each element of $%
\mathfrak{u}_{y}$ is invariant under both $G_{0}[y]$ and $\mathfrak{p}_{y}$.

Since a $G_{0}[y]-$invariant subspace of a polarized module need not be
invariant under the corresponding $\mathfrak{l}_{y}-$action, we define a%
\emph{\ morphism of polarized modules} to be a linear map that intertwines
both the $G_{0}[y]$ and $\mathfrak{p}_{y}-$actions. Thus the category of
polarized $G_{0}[y]-$modules is nothing but the category of
finite-dimensional $(\mathfrak{l}_{y},G_{0}[y])-$modules.

Let $Z(\mathfrak{l}_{y})$ denote the center of the enveloping algebra $U(%
\mathfrak{l}_{y})$. Since $\mathfrak{h}^{\ast }$ is the Cartan dual for $%
\mathfrak{l}_{y}$, the set of $\mathfrak{l}_{y}-$infinitesimal characters is
in natural correspondence with the quotient 
\begin{equation*}
\mathfrak{h}^{\ast }/W_{Y}
\end{equation*}%
where $W_{Y}$ is the Weyl group of $\mathfrak{h}^{\ast }$ in $\mathfrak{l}%
_{y}$. For $\lambda \in \mathfrak{h}^{\ast }$ we write $\sigma _{\lambda }$
to indicate the $\mathfrak{l}_{y}-$infinitesimal character corresponding to
the orbit $W_{y}\cdot \lambda $. A polarized $G_{0}[y]-$module $V$ is said
to have \emph{infinitesimal character }$\lambda \in \mathfrak{h}^{\ast }$ if 
$Z(\mathfrak{l}_{y})$ acts on $V$ by the character $\sigma _{\lambda }$.
Since $G_{0}$ is Harish-Chandra class, it follows that an irreducible
polarized $G_{0}[y]-$module has an infinitesimal character.

Given $y\in Y$ \ let $K[y]$ denote the stabilizer of $y$ in $K$. Then we can
introduce, making the obvious definitions, a category of polarized algebraic 
$K[y]-$modules, by adding the stipulation that the $K[y]-$action be
algebraic. Morphisms, as before, are linear maps that intertwine both the $%
K[y]$ and $\mathfrak{l}_{y}-$actions. We can also define and parametrize
infinitesimal characters as in the case of $G_{0}[y]$.

The following proposition can be deduced from the detailed description of
the stabilizers given in \cite{flag} via standard Lie theory considerations.

\begin{proposition}
Let $Y$ be a complex flag space for $G_{0}$ and suppose $y\in Y$ is special.
Then there exists a natural equivalence of categories between the category
of polarized $G_{0}[y]$-modules and the category of polarized algebraic $%
K[y] $-modules.
\end{proposition}

\bigskip

\noindent \textbf{The Standard Modules in Flag Spaces. }\ Suppose\ $y\in Y$
and let 
\begin{equation*}
\omega :G_{0}[y]\rightarrow GL(V)
\end{equation*}%
be an irreducible polarized representation. Let $S$ denote the $G_{0}-$orbit
of $y$. Then we have the corresponding homogeneous, analytic vector bundle

\begin{equation*}
\begin{array}{c}
\mathbb{V} \\ 
\downarrow \\ 
S%
\end{array}%
\end{equation*}%
with fiber $V$. The polarization allows us to define, in a canonical way,
the restricted holomorphic or \emph{polarized sections} of the analytic
vector bundle. In particular, let 
\begin{equation*}
\phi :G_{0}\rightarrow S\text{ \ \ be the projection }\phi (g)=g\cdot y\text{%
. }
\end{equation*}%
If $U\subseteq S$ is an open then a section of $\mathbb{V}$ over $U$ is a a
real analytic function 
\begin{equation*}
f:\phi ^{-1}(U)\rightarrow V\text{ \ such that \ }f(gp)=\omega (p^{-1})f(g)%
\text{ \ \ }\forall p\in G_{0}[y].
\end{equation*}%
The section is said to be \emph{polarized} if 
\begin{equation*}
\frac{d}{dt}\mid _{t=0}f(\exp (t\xi _{1})g)+i\frac{d}{dt}\mid _{t=0}f(\exp
(t\xi _{2})g)=-\omega (\xi _{1}+i\xi _{2})f(g)
\end{equation*}%
for all $\xi _{1}$, $\xi _{2}\in \mathfrak{g}_{0}$ such that $\xi _{1}+i\xi
_{2}\in \mathfrak{p}_{y}$.

Let $\mathcal{P}(y,V)$ denote the sheaf of polarized sections and let $%
\mathcal{O}_{Y}\mid _{S}$ be the sheaf of restricted holomorphic functions
on $S$. As a sheaf of $\mathcal{O}_{Y}\mid _{S}-$modules, $\mathcal{P}(y,V)$
is locally isomorphic to $\mathcal{O}_{Y}\mid _{S}\otimes V$ [cite]. The
left translation defines a $G_{0}$, and thus a $\mathfrak{g}-$action on $%
\mathcal{P}(y,V)$. Let $\lambda \in \mathfrak{h}^{\ast }$ be a parameter for
the $\mathfrak{l}_{y}-$infinitesimal character in $V$. Then the $\mathcal{O}%
_{Y}\mid _{S}$ and $\mathfrak{g}-$actions determine a $\pi _{\ast }(\mathcal{%
D}_{\lambda })\mid _{S}-$action. Put $\Theta =W\cdot \lambda $. Then the
compactly supported sheaf cohomology groups 
\begin{equation*}
H_{\text{c }}^{p}(S,\mathcal{P}(y,V))\text{ \ }p=0,1,2,3,\ldots
\end{equation*}%
are DNF $U_{\Theta }-$modules with a compatible $G_{0}-$action, provided
certain naturally defined topologies are Hausdorff \cite{analytic}.

Suppose $y\in S$ is special. Let $Q$ be the $K-$orbit of $y$ and let $q$ be
the codimension of the complex manifold $Q$ in $Y$. In general, one can show
the following. Although not difficult, the proof in \cite{preprint} depends
on some ideas which we will not otherwise use in this study.

\begin{proposition}
Maintain the above notations. \newline
\textbf{(a)} $H_{\text{c}}^{\text{p}}(S,\mathcal{P}(y,V))$ vanishes for $p<q$%
. \newline
\textbf{(b)} $H_{\text{c}}^{\text{n+q}}(S,\mathcal{P}(y,F))$ \ n$%
=0,1,2,\ldots $ is an admissible representation, naturally isomorphic to the
minimal globalization of its underlying Harish-Chandra module
\end{proposition}

Since $y$ is special, in a natural way $V$ is a polarized algebraic $K[y]-$%
module. Thus $V$ determines an algebraic vector bundle on the $K-$orbit $Q$.
The $\mathfrak{l}_{y}-$action in $V$, the translation by $K$, and the
natural $\mathcal{O}_{Q}-$action determine the action of a certain sheaf of
algebras, defined on $Q^{\text{alg}}$, on the corresponding sheaf of
algebraic sections \cite{chang}. Using a direct image construction \cite%
{chang}, modeled after the direct image for sheaves of TDOs modules \cite%
{HMSW}, one obtains a \emph{standard generalized Harish-Chandra sheaf }$%
\mathcal{I}(y,V)$ defined on the algebraic variety $Y^{\text{alg}}$. This
sheaf of $\mathcal{O}_{Y^{\text{alg}}}-$modules carries compatible actions
of $\mathfrak{g}$ and $K$. Indeed, $\mathcal{I}(y,V)$ is a sheaf of $\pi
_{\ast }(\mathcal{D}_{\lambda }^{\text{alg}})-$modules. One knows that the
corresponding sheaf cohomology groups 
\begin{equation*}
H^{\text{p}}(Y,\mathcal{I(}y,V))\text{ \ p}=0,1,2,\ldots
\end{equation*}%
are Harish-Chandra modules.

\bigskip

\noindent \textbf{Affinely Oriented Orbits.} \ A $K-$orbit $Q$ is called 
\emph{affinely embedded} if the inclusion 
\begin{equation*}
i:Q\rightarrow Y
\end{equation*}%
is an affine morphism. A $G_{0}-$orbit is called \emph{affinely oriented} if
its Matsuki dual is affinely embedded. Since the Matsuki dual of an open
orbit is closed \cite{matsuki}, it follows that all open $G_{0}-$orbits are
affinely oriented. It is known that all $K-$orbits in the full flag space
are affinely embedded, and more generally, if parabolic subalgebras in a $%
G_{0}-$orbit contain real Levi factors, then the orbit is affinely oriented 
\cite{chang2}.

By definition, \emph{a Levi orbit} is a $G_{0}-$orbit containing parabolic
subalgebras with real Levi factors. In the previous studies \cite{realizing}
and \cite{comparison2} only Levi orbits were considered. On the other hand,
it is not hard to define affinely embedded orbits which are not Levi. For
example, consider the natural action of the real special linear group $%
G_{0}=SL(n,\mathbb{R})$ on the complex projective space $Y=P^{n-1}(\mathbb{C}%
)$. If $n>2$, then there is a unique open $G_{0}-$orbit and this open orbit
is not Levi. In the last section of this paper we will consider a $G_{0}-$%
orbit which is not affinely oriented.

\bigskip

\noindent \textbf{Analytic Localization of Minimal Globalizations in the
Full Flag Space. \ }We conclude this section with the following theorem, due
to Hecht and Taylor, which characterizes the analytic localization to the
full flag space for the minimal globalization of a Harish-Chandra module
with regular infinitesimal character.

\begin{theorem}
Let $M$ be a Harish-Chandra module with regular infinitesimal character $%
\Theta $ and choose $\lambda \in \Theta $. Let $F_{\cdot }(M_{\text{min}})$
denote the Hochschild resolution for the minimal globalization of $M$ and
let 
\begin{equation*}
\Delta _{\lambda }(F_{\cdot }(M_{\text{min}}))
\end{equation*}
denote the corresponding analytic localization to the full flag space $X$.
Fix $x\in X$. Then we have the following. \newline
\textbf{(a)} Let $G_{0}[x]$ denote the stabilizer of $x$ in $G_{0}$. Then
the homology spaces of the complex 
\begin{equation*}
T_{x}\circ \Delta _{\lambda }(F_{\cdot }(M))
\end{equation*}
are finite-dimensional polarized $G_{0}[x]$-modules. \newline
\textbf{(b) }Let $S$ be the $G_{0}$-orbit of $x$ and let $h_{p}(\Delta
_{\lambda }(F_{\cdot }(M)))\mid _{S}$ denote the $p$-th homology of $\Delta
_{\lambda }(F_{\cdot }(M))$ restricted to $S$. Then $h_{p}(\Delta _{\lambda
}(F_{\cdot }(M)))\mid _{S}$ is the sheaf of polarized sections corresponding
to the polarized $G_{0}[x]$-module 
\begin{equation*}
h_{p}(T_{x}\circ \Delta _{\lambda }(F_{\cdot }(M))).
\end{equation*}
\end{theorem}

\noindent \textbf{Proof: }This result follows directly from\ Theorem 10.10,
Proposition 8.3 and Proposition 8.7 in \cite{analytic}. $\blacksquare $

\section{The Comparison Theorem}

In this section we generalize the Hecht-Taylor comparison theorem \cite%
{comparison1} to arbitrary orbits. In particular, suppose $y\in Y$ is
special, and let $\mathfrak{u}_{y}$ denote the nilpotent radical of the
corresponding parabolic subalgebra. We will establish the following theorem.

\begin{theorem}
\ Let $M$ be a Harish-Chandra module with regular infinitesimal character
and suppose $y$ is a special point. Assume that $M$\ has finite-dimensional $%
u_{y}-$homology groups and let $M_{\min }$ denote the minimal globalization
of $M$. Then, in a natural way, the Lie algebra homology groups 
\begin{equation*}
H_{p}(\mathfrak{u}_{y},M)\text{ \ and }H_{p}(\mathfrak{u}_{y},M_{\min })%
\text{ \ \ }p=0,1,2,\ldots
\end{equation*}%
are polarized $G_{0}[y]$-modules and the natural inclusion 
\begin{equation*}
M\rightarrow M_{\min }
\end{equation*}%
induces an isomorphism 
\begin{equation*}
H_{p}(\mathfrak{u}_{y},M)\cong H_{p}(\mathfrak{u}_{y},M_{\min })
\end{equation*}%
for each $p$.
\end{theorem}

\medskip

\noindent \textbf{Localization} \textbf{and} $\mathfrak{u}_{y}$-\textbf{%
homology}. \ Suppose $M$ is a $\mathfrak{u}_{y}-$module. By definition 
\begin{equation*}
H_{0}(\mathfrak{u}_{y},M)=\mathbb{C}\otimes _{\mathfrak{u}_{y}}M\text{.}
\end{equation*}%
When $M$ is a $\mathfrak{g}-$module then $H_{0}(\mathfrak{u}_{y},M)$ is a
module for the Levi quotient 
\begin{equation*}
\mathfrak{l}_{y}=\mathfrak{p}_{y}/\mathfrak{u}_{y}.
\end{equation*}%
The $\mathfrak{u}_{y}-$homology groups of $M$ are the derived functors of
the functor 
\begin{equation*}
M\mapsto H_{0}(\mathfrak{u}_{y},M)\text{.}
\end{equation*}%
Since $U(\mathfrak{g})$ is a free $U(\mathfrak{u}_{y})-$module, it follows
that a resolution of free $\mathfrak{g-}$modules can be used to compute the $%
\mathfrak{u}_{y}-$homology groups. In particular, if $M$ is a Harish-Chandra
module and if $F_{\cdot }(M)$ denotes the Hochschild resolution, then $K$
acts on $F_{\cdot }(M)$ by the tensor product of the adjoint action with the
action on $M$. This action is then homotopic to the left $\mathfrak{g}-$%
action. Thus one obtains a $K[y]-$action on the complex $H_{0}(\mathfrak{u}%
_{y},F_{\cdot }(M))$ and a corresponding algebraic $(\mathfrak{l}_{y},K[y])-$%
action on the $\mathfrak{u}_{y}-$homology groups. Similarly, there is a
continuous $G_{0}[y]-$action on the complex of DNF $\mathfrak{l}_{y}-$%
modules $H_{0}(\mathfrak{u}_{y},F_{\cdot }(M_{\text{min}}))$. Thus, since
these actions are homotopic, if the homology groups of $H_{0}(\mathfrak{u}%
_{y},F_{\cdot }(M_{\text{min}}))$ are finite-dimensional (and therefore
Hausdorff in the induced topologies), it follows that the homology spaces 
\begin{equation*}
H_{p}(\mathfrak{u}_{y},M_{\text{min}})\text{ \ }p=0,1,2,\ldots
\end{equation*}%
are polarized $G_{0}[y]-$modules. When $M$ is Harish-Chandra module with
infinitesimal character $\Theta $ then one can use the Hochschild resolution
with coefficients from $U_{\Theta }$ to compute the $\mathfrak{u}_{y}-$%
homology groups for $M$ and $M_{\text{min}}$, since $U_{\Theta }$ is a free $%
U(\mathfrak{u}_{y})-$module. The induced module structure on the homology
groups is independent of these two resolutions.

Let $Z(\mathfrak{l}_{y})$ denote the center of $U(\mathfrak{l}_{y})$ and
suppose $V$ is an $\mathfrak{l}_{y}-$module. For each $\lambda \in \mathfrak{%
h}^{\ast }$ we let $V_{\lambda }$ denote the corresponding $Z(\mathfrak{l}%
_{y})-$eigenspace in $V$. When $\Theta $ is a regular $\mathfrak{g}-$%
infinitesimal character and $M$ be a $U_{\Theta }-$module, then one knows
that 
\begin{equation*}
H_{p}(\mathfrak{u}_{y},M)=\bigoplus_{\lambda \in \Theta }H_{p}(\mathfrak{u}%
_{y},M)_{\lambda }.
\end{equation*}%
Indeed, letting $F_{\cdot }(M)$ denote the Hochschild resolution of $M$,
with coefficients from $U_{\Theta }$, one can deduce that the $p$-th
homology of the complex $H_{0}(\mathfrak{u}_{y},F_{\bullet }(M))_{\lambda }$
calculates the $\mathfrak{l}_{y}-$module $H_{p}(\mathfrak{u}_{y},M)_{\lambda
}$ \cite{comparison2}.

Thus, to establish the comparison theorem for a Harish-Chandra module with
regular infinitesimal character $\Theta $, it suffices to establish the
result for each of the spaces $H_{p}(\mathfrak{u}_{y},M)_{\lambda }$.

To calculate the modules $H_{p}(\mathfrak{u}_{y},M)_{\lambda }$, we use the
fact they can be identified with the derived functors of the geometric fiber
at $y$ of the corresponding localization to $Y$. We state this fact in the
following proposition. A proof can be found in \cite{comparison2}.

\begin{proposition}
Let $M$ be a $U_{\Theta }$-module with $\Theta $ regular and let $F_{\bullet
}(M)$ denote the corresponding Hochschild resolution of $M$. Choose $\lambda
\in \mathfrak{\Theta }$. Suppose $Y$ is a complex flag space and let $\Delta
_{\lambda }^{\text{alg}}$ denote the corresponding algebraic localization to 
$Y$. Then, for each $y\in Y$, there is a natural isomorphism of complexes of 
$\mathfrak{l}_{y}$-modules 
\begin{equation*}
T_{y}^{\text{alg}}\circ \Delta _{\lambda }^{\text{alg}}(F_{\bullet
}(M))\cong H_{0}(\mathfrak{u}_{y},F_{\bullet }(M))_{\lambda }.
\end{equation*}
\end{proposition}

\noindent \textbf{The Comparison Theorem} \ From the previous discussion,
the comparison theorem follows from the next result, which we prove in this
subsection.

\begin{theorem}
Let $M$ be a Harish-Chandra module with regular infinitesimal character $%
\Theta $. Suppose $y$ is a special point in a complex flag space $Y$ and let 
$\mathfrak{u}_{y}$ denote the nilradical of the corresponding parabolic
subalgebra $\mathfrak{p}_{y}$. Suppose $\lambda \in \Theta $ and assume that
each of the algebraic $(\mathfrak{l}_{y},K[y])$-modules 
\begin{equation*}
H_{p}(\mathfrak{u}_{y},M)_{\lambda }\text{ \ \ }p=0,1,2,\ldots
\end{equation*}%
is finite-dimensional. Then the natural inclusion 
\begin{equation*}
M\rightarrow M_{\min }
\end{equation*}%
induces an isomorphism 
\begin{equation*}
H_{p}(\mathfrak{u}_{y},M)_{\lambda }\cong H_{p}(\mathfrak{u}_{y},M_{\min
})_{\lambda }
\end{equation*}%
of polarized $G_{0}[y]$-modules, for each $p$.
\end{theorem}

\noindent \textbf{Proof: }\ Since the spaces $H_{p}(\mathfrak{u}%
_{y},M)_{\lambda }$ are finite-dimensional algebraic $(\mathfrak{l}%
_{y},K[y])-$modules and $y$ is special, these spaces are also polarized $%
G_{0}[y]-$modules. Indeed, in order to prove the theorem it suffices to show
that $H_{p}(\mathfrak{u}_{y},M)_{\lambda }$ and $H_{p}(\mathfrak{u}_{y},M_{%
\text{min}})_{\lambda }$ are isomorphic as $(\mathfrak{l}_{y},K_{0}[y])-$%
modules, where $K_{0}[y]$ is the stabilizer of $y$ in $K_{0}$.

Let $X$ be the full flag space and let $X_{y}$ be the fiber in $X$ over $y$.
Let $\Sigma _{Y}$ denote the root subspace of $\mathfrak{h}^{\ast }$
corresponding to Levi factors from $Y$ and let $W_{Y}$ denote the associated
Weyl group. Put 
\begin{equation*}
\Sigma _{Y}^{+}=\Sigma _{Y}\cap \Sigma ^{+}.
\end{equation*}%
We say that $\lambda $ is \emph{antidominant for the fiber }if \ 
\begin{equation*}
\overset{\vee }{\alpha }(\lambda )\notin \left\{ 1,2,3,\ldots \right\} \text{
\ \ for each }\alpha \in \Sigma _{Y}^{+}\text{.}
\end{equation*}%
Since there exists $w\in W_{Y}$ such that $w\lambda $ is antidominant for
the fiber and since $w\lambda $ and $\lambda $ parameterize the same $%
\mathfrak{l}_{y}-$infinitesimal character, we may assume that $\lambda $ is
antidominant for the fiber. Suppose $x\in X_{y}$ and let 
\begin{equation*}
i:X_{y}\rightarrow X
\end{equation*}%
denote the inclusion. Reintroducing the notations established in Section 2,
we now prove the following lemma.

\medskip

\noindent \textbf{Lemma }\emph{\ Maintaining the assumptions of Theorem 4.3,
let} 
\begin{equation*}
h_{p}(i^{\ast }\circ \Delta _{\lambda ,X}(F_{\cdot }(M))
\end{equation*}%
\emph{denote the} $p$\emph{-th homology of the complex} $i^{\ast }\circ
\Delta _{\lambda ,X}(F_{\cdot }(M)$. \emph{Then} $h_{p}(i^{\ast }\circ
\Delta _{\lambda ,X}(F_{\cdot }(M))$\emph{\ is the sheaf of holomorphic
section of a }$K[y]$-\emph{equivariant} \emph{holomorphic vector bundle over}
$X_{y}$.

\medskip

\noindent \textbf{Proof of Lemma: \ }By Proposition 2.2 and Proposition 4.2,
it follows from the given assumptions that the homology groups of the
complex 
\begin{equation*}
\Gamma (X_{y},i_{\text{alg}}^{\ast }\circ \Delta _{\lambda ,X}^{\text{alg}%
}(F_{\cdot }(M)))
\end{equation*}%
are finite-dimensional algebraic $(\mathfrak{l}_{y},K_{y})-$modules. We
claim that this implies that the homologies of the complex 
\begin{equation*}
i_{\text{alg}}^{\ast }\circ \Delta _{\lambda ,X}^{\text{alg}}(F_{\cdot }(M))
\end{equation*}%
are the sheaves of sections for $K[y]-$equivariant algebraic vector bundles
defined over the algebraic variety $X_{y}^{\text{alg}}$. In particular, one
knows that $X_{y}^{\text{alg}}$ is the full flag space for the Levi quotient 
$\mathfrak{l}_{y}$ and that the homology groups of the previous complex%
\textbf{\ }are sheaves of modules for a twisted sheaf of differential
operators $\mathcal{D}_{\lambda ,X_{y}^{\text{alg}}}^{\text{alg}}$ defined
on $X_{y}^{\text{alg}}$. Since the parameter $\lambda $ is antidominant with
respect to $\mathfrak{l}_{y}$, it follows that the global sections define an
exact functor on the category of quasicoherent $\mathcal{D}_{\lambda ,X_{y}^{%
\text{alg}}}^{\text{alg}}-$modules. Thus, for each $p=0,1,2,\ldots $, there
are natural isomorphisms 
\begin{equation*}
h_{p}\left( \Gamma (X_{y},i_{\text{alg}}^{\ast }\circ \Delta _{\lambda ,X}^{%
\text{alg}}(F_{\cdot }(M)))\right) \cong \Gamma (X_{y},h_{p}(i_{\text{alg}%
}^{\ast }\circ \Delta _{\lambda ,X}^{\text{alg}}(F_{\cdot }(M))))
\end{equation*}%
of finite-dimensional algebraic $(\mathfrak{l}_{y},K[y])-$modules, where $%
h_{p}(\cdot )$ denotes the $p$-th homology group of the given complex.
Therefore, our claim follows, since the only quasicoherent sheaves of $%
\mathcal{D}_{\lambda ,X_{y}^{\text{alg}}}^{\text{alg}}-$modules with
finite-dimensional global sections are finite-rank locally free sheaves of $%
\mathcal{O}_{X_{y}^{\text{alg}}}-$modules.

Now let $j:X_{y}\rightarrow X_{y}^{\text{alg}}$ indicate the identity and
let 
\begin{equation*}
\epsilon \left( \cdot \right) =\mathcal{O}_{X_{y}}\otimes _{j^{-1}(\mathcal{O%
}_{X_{y}^{\text{alg}}})}j^{-1}(\cdot )
\end{equation*}%
denote Serre's GAGA functor \cite{serre}. We claim that 
\begin{equation*}
\epsilon \circ h_{p}(i_{\text{alg}}^{\ast }\circ \Delta _{\lambda ,X}^{\text{%
alg}}(F_{\cdot }(M)))\cong h_{p}\left( i^{\ast }\circ \Delta _{\lambda
,X}(F_{\cdot }(M))\right) \text{.}
\end{equation*}%
Indeed, the claim follows, since there is a natural isomorphism of complexes
of sheaves of $(\mathfrak{l}_{y},K[y])-$modules 
\begin{equation*}
\epsilon \circ i_{\text{alg}}^{\ast }\circ \Delta _{\lambda ,X}^{\text{alg}%
}(F_{\cdot }(M))\cong i^{\ast }\circ \Delta _{\lambda ,X}(F_{\cdot }(M)).
\end{equation*}%
and since the functor $\epsilon $ is exact on the category of quasicoherent $%
\mathcal{O}_{X_{y}^{\text{alg}}}-$modules. This proves the lemma. $%
\blacksquare $

\medskip

\noindent We now continue with the proof of Theorem 4.3, by establishing the
following lemma.

\medskip

\noindent \textbf{Lemma \ }\emph{Use the given notations and maintain the
assumptions of Theorem 4.3. Then the natural morphism } 
\begin{equation*}
i^{\ast }\circ \Delta _{\lambda ,X}(F_{\cdot }(M))\rightarrow i^{\ast }\circ
\Delta _{\lambda ,X}(F_{\cdot }(M_{\text{min}}))
\end{equation*}%
\emph{of complexes of sheaves of} $(\mathfrak{l}_{y},K_{0}[y])$\emph{%
-modules, induces an isomorphism on the level of homology groups. }

\medskip

\noindent \textbf{Proof of Lemma: }It follows from Theorem 3.3, that for
each $x\in X$, the stalks of the of the homology sheaves $h_{p}\left( \Delta
_{\lambda ,X}(F_{\cdot }(M_{\text{min}})\right) $ are locally free, finite
rank $\mathcal{O}_{X}\mid _{x}-$modules. Therefore, for each $x\in X_{y}$,
the homology sheaves 
\begin{equation*}
h_{p}\left( i^{\ast }\circ \Delta _{\lambda ,X}(F_{\cdot }(M_{\text{min}%
})\right)
\end{equation*}%
are locally free, finite rank $\mathcal{O}_{X_{y}}\mid _{x}-$modules. We now
apply the comparison theorem of Hecht and Taylor \cite{comparison1} to
deduce the desired isomorphism. Specifically: for $x\in X_{y}$, let $%
T_{x,X_{y}}$ denote the functor that takes the geometric fiber at $x$ with
respect to sheaves of $\mathcal{O}_{X_{y}}-$modules. Then the Hecht-Taylor
result implies that the natural morphism 
\begin{equation*}
T_{x,X_{y}}\circ i^{\ast }\circ \Delta _{\lambda ,X}(F_{\cdot
}(M))\rightarrow T_{x,X_{y}}\circ i^{\ast }\circ \Delta _{\lambda
,X}(F_{\cdot }(M_{\text{min}}))
\end{equation*}%
induces an isomorphism on homology groups, when $x\in X_{y}$ is special.
Thus for each special point $x\in X_{y}$ and for each whole number $p$, we
have a natural isomorphism 
\begin{equation*}
T_{x,X_{y}}\circ h_{p}(i^{\ast }\circ \Delta _{\lambda ,X}(F_{\cdot
}(M)))\cong T_{x,X_{y}}\circ h_{p}\left( i^{\ast }\circ \Delta _{\lambda
,X}(F_{\cdot }(M_{\text{min}}))\right) \text{.}
\end{equation*}%
Therefore the lemma follows, since there is a special point in each $%
G_{0}[y]-$orbit on $X_{y}$ \cite{matsuki}. $\blacksquare $

\medskip

\noindent We can now conclude that the natural morphism 
\begin{equation*}
\Gamma (X_{y},i^{\ast }\circ \Delta _{\lambda ,X}(F_{\cdot }(M)))\rightarrow
\Gamma (X_{y},i^{\ast }\circ \Delta _{\lambda ,X}(F_{\cdot }(M_{\text{min }%
})))
\end{equation*}%
induces an isomorphism on the level of homology groups as follows. Therefore
the proof of Theorem 4.3 follows immediately by an application of
Proposition 2.2, Proposition 2.1 and Proposition 4.2. $\ \blacksquare $

\section{Geometric Realization of Representations}

In this section we consider the relation of the comparison theorem to the
geometric realization of representations in complex flag spaces. In
particular, suppose $y\in Y$ \ is a special point and let $V$ be an
irreducible polarized $G_{0}[y]-$module. Let $\mathcal{I}(y,V)$ denote the
corresponding generalized standard Harish-Chandra sheaf defined on $Y^{\text{%
alg}}$. If \ $\lambda \in \mathfrak{h}^{\ast }$ is a parameter for the $%
\mathfrak{l}_{y}-$infinitesimal character in $V$ then the sheaf cohomologies 
\begin{equation*}
H^{p}(Y,\mathcal{I}(y,F))\text{ \ \ }p=0,1,2,\ldots
\end{equation*}%
are Harish-Chandra modules with $\mathfrak{g}-$infinitesimal character $%
\Theta =W\cdot \lambda $. Put 
\begin{equation*}
M=\Gamma (Y,\mathcal{I}(y,F))
\end{equation*}%
and let $M_{\text{min}}$ denote the minimal globalization of $M$. We are
interested in finding a geometric realization for $M_{\text{min}}$ in $Y$.
One obvious candidate is the analytic localization of $M_{\text{min}}$ to $Y$%
. In fact, suppose $\mu \in \Theta $ and, using our previously established
notation, let 
\begin{equation*}
\Delta _{\mu }(F_{\cdot }(M_{\text{min}}))
\end{equation*}%
denote the analytic localization of $M_{\min }$ to $Y$. It follows from the
Beilinson-Bernstein result that the sheaves $\Delta _{\mu }(F_{p}(M_{\text{%
min}}))$ $p=0,1,2,\ldots $ are acyclic for the functor of global sections
and that there is a natural isomorphism of complexes 
\begin{equation*}
\Gamma (Y,\Delta _{\mu }(F_{\cdot }(M_{\text{min}}))\cong F_{\cdot }(M_{%
\text{min}}).
\end{equation*}%
Thus the complex $\Delta _{\mu }(F_{\cdot }(M_{\text{min}}))$ has vanishing
hypercohomology in all degrees except zero, where we reobtain the module $M_{%
\text{min}}$. Indeed, when the infinitesimal character $\Theta $ is regular,
one obtains the following uniqueness for this geometric realization of $M_{%
\text{min}}$. Let $\mathcal{F}_{\cdot }$ be a complex of sheaves of DNF $\pi
_{\ast }(\mathcal{D}_{\mu }).-$modules, with bounded homology, whose
hypercohomology realizes the module $M_{\text{min}}$, then there are natural
isomorphisms in homology 
\begin{equation*}
h_{p}(\Delta _{\mu }(F_{\cdot }(M_{\text{min}})))\cong h_{p}(\mathcal{F}%
_{\cdot })
\end{equation*}%
for each $p$. In the case of the full flag space, this uniqueness follows
from an equivalence of derived categories shown in \cite{analytic}, and the
general case is not hard to deduce from this.

Thus we would like to understand the complex $\Delta _{\mu }(F_{\cdot }(M_{%
\text{min}}))$. It turns out that the structure of the analytic localization
is completely determined by the corresponding geometric fibers. In
particular, we have the following result \cite{realizing}.

\begin{proposition}
Let $W$ be a minimal globalization with $\mathfrak{g}-$infinitesimal
character $\Theta $ and choose $\lambda \in \Theta $. Suppose $Y$ is a
complex flag space for $G_{0}$ and, using the previously established
notation, let $\Delta _{\lambda }(F_{\cdot }(W))$ denote the analytic
localization of $W$ to $Y$. Choose $y\in Y$ and let $S=G_{0}\cdot y$. Assume
that each of the homology groups 
\begin{equation*}
h_{p}(T_{y}\circ \Delta _{\lambda }(F_{\cdot }(W)))\text{ \ }p=0,1,2,\ldots
\end{equation*}%
is finite-dimensional. Let $\mathcal{P}(y,h_{p}(T_{y}\circ \Delta _{\lambda
}(F_{\cdot }(W))))$ denote the sheaf of polarized sections for the polarized
homogeneous vector bundle on $S$ determined by $h_{p}(T_{y}\circ \Delta
_{\lambda }(F_{\cdot }(W)))$. Then there is a natural isomorphism 
\begin{equation*}
h_{p}(\Delta _{\lambda }(F_{\cdot }(W)))\mid _{S}\cong \mathcal{P}%
(y,h_{p}(T_{y}\circ \Delta _{\lambda }(F_{\cdot }(W))))\text{ }
\end{equation*}%
of $G_{0}$-equivariant DNF sheaves of $\mathfrak{g}-$modules.
\end{proposition}

\noindent \textbf{Base Change \ }The previous proposition, in conjunction
with the comparison theorem, can be used to deduce information about the
geometric realization for the minimal globalization of a generalized
standard Beilinson-Bernstein module. There is a third ingredient we will
need for our analysis: the so-called base change formula \cite{realizing},
as applied to the derived geometric fibers of the Harish-Chandra sheaf $%
\mathcal{I}(y,V)$. In particular, let $D^{b}(\pi _{\ast }(\mathcal{D}%
_{\lambda }^{\text{alg}}))$ (respectively $D^{b}(U_{\lambda }(\mathfrak{l}%
_{y})$) denote the derived category of bounded complexes of quasi-coherent $%
\pi _{\ast }(\mathcal{D}_{\lambda }^{\text{alg}})-$modules (respectively the
derived category of bounded complexes of $\mathfrak{l}_{y}-$modules with
infinitesimal character $\lambda $). Suppose $z\in Y$. Then, in a natural
way, the geometric fiber at $z$ determines a derived functor 
\begin{equation*}
LT_{z}:D^{b}(\pi _{\ast }(\mathcal{D}_{\lambda }^{\text{alg}}))\rightarrow
D^{b}(U_{\lambda }(\mathfrak{l}_{y})).
\end{equation*}%
Let $Q$ be the $K-$orbit of $y$ and let $\overline{Q}$ be the Zariski
closure of $Q$ in $Y$. Put $\partial Q=\overline{Q}-Q$ and $U=Y-\partial Q$.
Thus $U$ is Zariski open. Let $q$ denote the codimension of $Q$ in $Y$ and
let $V[q]$ denote the complex of $\mathfrak{l}_{y}-$modules which is zero
except in homology degree $q$, where one obtains the module $V$. We also
identify the sheaf $\mathcal{I}(y,V)$ with the complex in $D^{b}(\pi _{\ast
}(\mathcal{D}_{\lambda }^{\text{alg}}))$ which is zero in all degrees except
degree zero where we obtain $\mathcal{I}(y,V)$. Then, at least for $z\in U$,
the complex $LT_{z}(\mathcal{I}(y,V))$ is simple to understand. We summarize
in the following proposition.

\begin{proposition}
Maintain the previously introduced notations. Then we have the following
isomorphisms in $D^{b}(U_{\lambda }(\mathfrak{l}_{y}))$. \newline
(\textbf{a}) For $z\in U-Q$ 
\begin{equation*}
LT_{z}(\mathcal{I}(y,V))\cong 0.
\end{equation*}%
\newline
(\textbf{b}) $LT_{y}(\mathcal{I}(y,V))\cong V[q]$.
\end{proposition}

\textbf{Proof: }The result follows from the construction of $\mathcal{I}%
(y,F) $ and the base change formula, which holds for the generalized direct
image, as in the case of the direct image for $\mathcal{D}-$modules \cite%
{Borel}. $\blacksquare $

\medskip

For $z\in \partial Q$, the structure the complex $LT_{z}(\mathcal{I}(y,V))$
is more complicated, at least when $Q$ is not affinely embedded in $Y$. In
particular, let 
\begin{equation*}
i:U\rightarrow Y
\end{equation*}%
denote the inclusion. We let $\pi _{\ast }(\mathcal{D}_{\lambda }^{\text{alg}%
})\mid _{U}$ be the sheaf of algebras $\pi _{\ast }(\mathcal{D}_{\lambda }^{%
\text{alg}})$ restricted to $U$ and let $D^{b}(\pi _{\ast }(\mathcal{D}%
_{\lambda }^{\text{alg}})\mid _{U})$ denote the derived category of bounded
complexes of quasi-coherent $\pi _{\ast }(\mathcal{D}_{\lambda }^{\text{alg}%
})\mid _{U}$-modules. Then the direct image in the category of sheaves
induces a derived functor 
\begin{equation*}
Ri_{\ast }:D^{b}(\pi _{\ast }(\mathcal{D}_{\lambda }^{\text{alg}})\mid
_{U})\rightarrow D^{b}(\pi _{\ast }(\mathcal{D}_{\lambda }^{\text{alg}})).
\end{equation*}%
We have the following.

\begin{proposition}
Maintain the previously introduced notations. \newline
(a) Suppose $z\in \partial Q$. Then 
\begin{equation*}
LT_{z}\circ Ri_{\ast }(\mathcal{I}(y,V)\mid _{U})\cong 0
\end{equation*}%
in the category $D^{b}(U_{\lambda }(\mathfrak{l}_{y}))$. \newline
(b) If $Q$ is affinely embedded in $Y$ then 
\begin{equation*}
Ri_{\ast }(\mathcal{I}(y,V)\mid _{U})\cong \mathcal{I}(y,V).
\end{equation*}%
In particular 
\begin{equation*}
LT_{z}(\mathcal{I}(y,V))\cong 0
\end{equation*}%
for $z\in \partial Q$.
\end{proposition}

\textbf{Proof:} \ Once again, the first claim (a) is an application of the
base change formula for the generalized direct image, applied to the sheaf $%
\mathcal{I}(y,V)$. On the other hand, the second claim (b) is another
standard result for the direct image functor \cite{Borel}, which also
applies to the direct image in the category of generalized $\mathcal{D}-$%
modules. $\blacksquare $

\medskip

\noindent \textbf{Geometric Realization for the Minimal Globalization of a
Standard Module }\ Suppose 
\begin{equation*}
M=\Gamma (Y,\mathcal{I}(y,V))
\end{equation*}%
is a standard Harish-Chandra module, where $y\in Y$ is special and $V$ is an
irreducible polarized $G_{0}[y]-$module. Let $\mathcal{P}(y,V)$ be the sheaf
of polarized sections for the corresponding $G_{0}$ homogeneous polarized
vector bundle and let $M_{\text{min}}$ denote the minimal globalization of $%
M $. We are now ready to deduce the following result, which generalizes the
result for the full flag space.

\begin{theorem}
Maintain the previous assumptions and notations. Let $Q$ denote the $K-$%
orbit of the special point $y$ and let $q$ denote the codimension of $Q$ in $%
Y$. Assume the $\mathfrak{l}_{y}$-infinitesimal character in $V$ is regular
and antidominant for $Y$ and let $\lambda \in \mathfrak{h}^{\ast }$ be a
corresponding parameter. \newline
(\textbf{a}) Suppose $S$ is the $G_{0}$-orbit of $y$. Then 
\begin{equation*}
h_{p}\left( \Delta _{\lambda }(F_{\cdot }(M_{\text{min}}))\right) \mid
_{S}\cong 
\begin{array}{c}
0\;\;\text{for }p\neq q \\ 
\mathcal{P}(y,F)\text{ \ for }p=q%
\end{array}%
\text{ .}
\end{equation*}%
\newline
\textbf{(b)} Suppose $S$ is affinely oriented. Then 
\begin{equation*}
h_{p}\left( \Delta _{\lambda }(F_{\cdot }(M_{\text{min}}))\right) \cong 
\begin{array}{c}
0\;\;\text{for }p\neq q \\ 
\mathcal{P}(y,F)^{Y}\text{ \ }p=q%
\end{array}%
\text{ }
\end{equation*}%
where $\mathcal{P}(y,F)^{Y}$ denotes the extension by zero of $\mathcal{P}%
(y,F)$ to $Y$.
\end{theorem}

\textbf{Proof:} Since $\lambda $ is regular and antidominant for $Y$, it
follows from the Beilnson-Bernstein equivalence of categories that 
\begin{equation*}
h_{p}\left( \Delta _{\lambda }^{\text{alg}}(F_{\cdot }(M))\right) \cong 
\begin{array}{c}
0\text{ \ for }p\neq 0 \\ 
\mathcal{I}(y,V)\text{ \ \ for \ }p=0%
\end{array}%
.
\end{equation*}%
Thus, for $z\in Y$, the homologies of the complex $T_{z}\circ \Delta
_{\lambda }^{\text{alg}}(F_{\cdot }(M))$ are isomorphic to the homologies of 
$LT_{z}(\mathcal{I}(y,V))$. Via the comparison theorem, the homology groups
of $T_{z}\circ \Delta _{\lambda }^{\text{alg}}(F_{\cdot }(M))$ coincide with
the homology groups of \ $T_{z}\circ \Delta _{\lambda }(F_{\cdot }(M_{\text{%
min}}))$ when these homology groups are finite dimensional. Thus the first
part of the theorem follows by an application Proposition 5.2 together with
the Proposition 5.1, and the second part follows easily using Proposition
5.3. $\blacksquare $

\begin{theorem}
Suppose $y\in Y$ is special and that the $K-$orbit $Q$ of $y$ is affinely
embedded in $Y.$ Suppose $V$ is an irreducible polarized $G_{0}[y]$-module.
Let $\mathcal{I}(y,V)$ indicate the corresponding standard Harish-Chandra
sheaf on $Y^{\text{alg}}$ and let $\mathcal{P}(y,V)$ denote the
corresponding sheaf of polarized sections on the $G_{0}$-orbit $S$ of $y$.
Suppose $q$ is the codimension of $Q$ in $Y$. Then the compactly supported
cohomology group $H_{\text{c}}^{p}(S,\mathcal{P}(y,V))$ vanishes for $p<q$
and for each $n\geq 0$, $H_{\text{c}}^{q+n}(S,\mathcal{P}(y,V))$ is
naturally isomorphic to the minimal globalization of $H^{n}(Y,\mathcal{I}%
(y,V))$.
\end{theorem}

\textbf{Proof:} \ When $V$ has an $\mathfrak{l}_{y}-$infinitesimal character
that is regular and antidominant for $Y$, then the corollary\ follows
immediately from the previous theorem. The general case follows by a
tensoring argument, as in \cite{realizing}. $\blacksquare $\noindent 

\section{The $SU(n,1)$-action in Complex Projective Space}

In this section we give an example to analyze the situation when the $G_{0}-$%
orbit is not affinely oriented. In particular, we show that Theorem 5.5
fails to hold. Put 
\begin{equation*}
J=\sum_{j=1}^{n}E_{jj}-E_{n+1n+1}
\end{equation*}%
where $E_{jk}$ are the standard basis for the $(n+1)\times (n+1)$ matrices
and suppose $G$ is the complex special linear group $SL(n+1,\mathbb{C})$.
Define

\begin{equation*}
\gamma (A)=(\overline{A}^{t})^{-1}\text{ \ and \ }\tau (A)=J\gamma (A)J
\end{equation*}%
for $A\in G$. Thus $\tau $ is a conjugation of $G$ and $\gamma $ is a
compact conjugation commuting with $\tau $. By definition, the fixed point
set of $\tau $ is the group 
\begin{equation*}
G_{0}=SU(n,1).
\end{equation*}%
The corresponding $\gamma -$invariant maximal compact subgroup of $G_{0}$ is 
\begin{equation*}
K_{0}=SU(n+1)\cap G_{0}\text{.}
\end{equation*}%
The complexification $K$ of $K_{0}$ is naturally isomorphic to the fixed
point in $G$ of the involution 
\begin{equation*}
\theta (A)=JAJ.
\end{equation*}%
Thus the elements of $K$ are the matrices of the form 
\begin{equation*}
\left( 
\begin{array}{cccc}
&  &  & 0 \\ 
& A &  & \vdots \\ 
&  &  & 0 \\ 
0 & \cdots & 0 & \left( \text{det}A\right) ^{-1}%
\end{array}%
\right)
\end{equation*}%
where $A\in GL(n,\mathbb{C})$.

We calculate the $K-$orbits on the complex flag space $Y=P^{n}(\mathbb{C})$.
For $(z_{1},\ldots ,z_{n+1})\in \mathbb{C}^{n+1}$ let 
\begin{equation*}
\left[ 
\begin{array}{c}
z_{1} \\ 
\vdots \\ 
z_{n+1}%
\end{array}%
\right]
\end{equation*}%
denote the corresponding point in $P^{n}(\mathbb{C})$. Let $U\subseteq P^{n}(%
\mathbb{C})$ be the $K-$invariant affine open set defined by $z_{n+1}\neq 0$%
. Thus $U$ contains two $K-$orbits: one consisting of a fixed point 
\begin{equation*}
Q_{\text{fp}}=\left\{ \left[ 
\begin{array}{c}
0 \\ 
\vdots \\ 
0 \\ 
1%
\end{array}%
\right] \right\}
\end{equation*}%
and the other being the open $K-$orbit: 
\begin{equation*}
Q_{\text{o}}=U-Q_{\text{fp}}\text{.}
\end{equation*}%
The complement of $U$: 
\begin{equation*}
Q_{\text{c}}=P^{n}(\mathbb{C})-U
\end{equation*}%
is a closed $K-$orbit of dimension $n-1$.

Matsuki duality now determines the $G_{0}-$orbits on $Y$. In particular, let 
$S_{\text{fp}}$, $S_{\text{o}}$ and $S_{\text{c}}$ denote the dual orbits to 
$Q_{\text{fp}}$, $Q_{\text{o}}$ and $Q_{\text{c}}$, respectively. Thus $S_{%
\text{fp}}$ and $S_{\text{c}}$ are open orbits while $S_{\text{o}}$ is
closed in $Y$. Observe that $S_{\text{fp}}$ and $S_{\text{c}}$ are affinely
oriented while $S_{\text{o}}$ is not when $n>1$. Using Iwasawa decomposition
for $G_{0}$, one sees that $K_{0}$ acts transitively on $S_{\text{o}}$. In
particular, $S_{\text{o}}\subseteq Q_{\text{o}}$ and each point in $S_{\text{%
o}}$ is special.

Let $\mathcal{O}_{Y}$ denote the sheaf of holomorphic functions on $Y$ and
let $\mathcal{O}_{Y}\mid _{S_{\text{o}}}$ denote the restriction of $%
\mathcal{O}_{Y}$ to $S_{\text{o}}$. We also introduce sheaf $\mathcal{O}_{Y^{%
\text{alg}}}$ of regular functions on the algebraic variety $Y^{\text{alg}}$
and let 
\begin{equation*}
i:Q_{\text{o}}^{\text{alg}}\rightarrow Y^{\text{alg}}
\end{equation*}%
denote the inclusion. Choose a point $y\in S_{\text{o}}$ and let $\mathbb{C}$
denote trivial one-dimensional polarized $G_{0}[y]$-module. Then the
corresponding sheaf $\mathcal{P}(y,\mathbb{C)}$ of polarized sections is the 
$G_{0}-$equivariant sheaf $\mathcal{O}_{Y}\mid _{S_{\text{o}}}$and the
corresponding Harish-Chandra sheaf $\mathcal{I(}y,\mathbb{C})$ is the $K-$%
equivariant sheaf of $\mathfrak{g}-$modules $i_{\ast }(\mathcal{O}_{Y^{\text{%
alg}}}\mid _{Q_{\text{o}}^{\text{alg}}})$ where $i_{\ast }$ denotes the
direct image in the category of sheaves.

Let $\mathfrak{u}_{y}$ denote the nilradical of the parabolic subalgebra $%
\mathfrak{p}_{y}$ and let 
\begin{equation*}
\mathfrak{l}_{y}=\mathfrak{p}_{y}/\mathfrak{u}_{y}
\end{equation*}%
denote the corresponding Levi quotient. Then the $\mathfrak{l}_{y}-$%
infinitesimal character for the trivial module $\mathbb{C}$ is parametrized
by $-\rho $, where $\rho $ is one half the sum of the positive roots in $%
\mathfrak{h}^{\ast }$. Since $-\rho $ is regular and antidominant, it
follows that 
\begin{equation*}
H^{p}(Y,\mathcal{I(}y,\mathbb{C}))=0\text{ \ for }p>0\text{ \ and \ }\Gamma
(Y,\mathcal{I(}y,\mathbb{C}))\neq 0.
\end{equation*}%
By a direct calculation, it is not hard to show that the set of $K_{0}-$%
finite vectors in 
\begin{equation*}
\Gamma (S_{\text{o}},\mathcal{O}_{Y}\mid _{S_{\text{o}}})=\Gamma (S_{\text{o}%
},\mathcal{P}(y,\mathbb{C))}
\end{equation*}%
is naturally isomorphic to 
\begin{equation*}
\Gamma (Q_{\text{o}},\mathcal{O}_{Y^{\text{alg}}}\mid _{Q_{\text{o}}^{\text{%
alg}}})=\Gamma (Y,\mathcal{I(}y,\mathbb{C}))
\end{equation*}%
although we shall give a different reason for this below. On the other hand,
since the codimension of $Q_{\text{o}}$ in $Y$ is zero and since $S_{\text{o}%
}$ is compact, if the orbit $Q_{\text{o}}$ were affinely imbedded, it would
follow from the work in the last section of this paper that 
\begin{equation*}
H^{p}(S_{\text{o}},\mathcal{P(}y,\mathbb{C}))=0\text{ \ for }p>0\text{.\ }
\end{equation*}%
We now show that this vanishing does not occur when $n>1$.

If $\mathcal{F}$ is a sheaf defined on a locally closed subset of $Y$, we
let $\mathcal{F}^{Y}$ denote the extension by zero of $\mathcal{F}$ to $Y$.
To calculate the higher sheaf cohomologies of $\mathcal{O}_{Y}\mid _{S_{%
\text{o}}}$ we consider the following short exact sequence of $G_{0}-$%
equivariant sheaves on $Y$: 
\begin{equation*}
0\rightarrow \left( \mathcal{O}_{Y}\mid _{S_{\text{c}}\cup S_{\text{fp}%
}}\right) ^{Y}\rightarrow \mathcal{O}_{Y}\rightarrow \left( \mathcal{O}%
_{Y}\mid _{S_{\text{o}}}\right) ^{Y}\rightarrow 0.
\end{equation*}%
We compute the resulting long exact sequence in sheaf cohomology. Since $S_{%
\text{c}}$ and $S_{\text{fp}}$ are open orbits, a standard sheaf cohomology
argument shows that 
\begin{equation*}
H^{p}(Y,\left( \mathcal{O}_{Y}\mid _{S_{\text{c}}\cup S_{\text{fp}}}\right)
^{Y})\cong H_{\text{c}}^{p}(S_{\text{c}},\mathcal{O}_{S_{\text{c}}})\oplus
H_{\text{c}}^{p}(S_{\text{fp}},\mathcal{O}_{S_{\text{fp}}})
\end{equation*}%
where $\mathcal{O}_{S_{\text{c}}}$ and $\mathcal{O}_{S_{\text{fp}}}$ denote
the sheaves of holomorphic functions on $S_{\text{c}}$ and $S_{\text{fp}}$,
respectively. Since the codimension of $Q_{\text{c}}$ is one and the
codimension of $Q_{\text{fp}}$ is $n$ it follows from Theorem 5.5 that 
\begin{equation*}
H_{\text{c}}^{p}(S_{\text{c}},\mathcal{O}_{S_{\text{c}}})=0\text{ for }p\neq
1\text{ \ and \ }H_{\text{c}}^{p}(S_{\text{fp}},\mathcal{O}_{S_{\text{fp}%
}})=0\text{ \ for }p\neq n
\end{equation*}%
Indeed, via Kashiwara's equivalence of categories for the direct image
functor \cite{chang} one deduces that $H_{\text{c}}^{1}(S_{\text{c}},%
\mathcal{O}_{S_{\text{c}}})$ and $H_{\text{c}}^{n}(S_{\text{fp}},\mathcal{O}%
_{S_{\text{fp}}})$ are irreducible minimal globalizations. In particular,
each of these last two cohomologies are nonzero. On the other hand, it is
well known that sheaf cohomology for $\mathcal{O}_{Y}$ vanishes in positive
degree. Thus, for $n>1$, we obtain the short exact sequence 
\begin{equation*}
0\rightarrow \mathbb{C\rightarrow }\Gamma (S_{\text{o}},\mathcal{P}(y,%
\mathbb{C)})\rightarrow H_{\text{c}}^{1}(S_{\text{c}},\mathcal{O}_{S_{\text{c%
}}})\rightarrow 0.
\end{equation*}%
For positive $p$, it follows that $H^{p}(S_{\text{o}},\mathcal{P}(y,\mathbb{%
C)})$ is zero except when $p=n$, in which case we obtain the isomorphism: 
\begin{equation*}
H^{n-1}(S_{\text{o}},\mathcal{P}(y,\mathbb{C)})\cong H_{\text{c}}^{n}(S_{%
\text{fp}},\mathcal{O}_{S_{\text{fp}}})
\end{equation*}%
which contradicts Theorem 5.5.

We continue analyzing this example using ideas developed in our study. Put 
\begin{equation*}
M=\Gamma (Y,\mathcal{I}(y,\mathbb{C))}
\end{equation*}%
and let $M_{\text{min}}$ denote the minimal globalization of $M$. We
calculate the analytic localization of $M_{\text{min}}$ to $Y$ and use this
information to deduce that 
\begin{equation*}
M_{\text{min}}\cong \Gamma (S,\mathcal{P}(y,\mathbb{C})).
\end{equation*}%
Let 
\begin{equation*}
\Delta _{-\rho }(F_{\cdot }(M_{\text{min}}))
\end{equation*}%
denote the corresponding analytic localization of $M_{\text{min}}$ to $Y$.
By the comparison theorem, for $z\in Y$ special, the morphism of complexes 
\begin{equation*}
T_{z}^{\text{alg}}\circ \Delta _{-\rho }^{\text{alg}}(F_{\cdot
}(M))\rightarrow T_{z}\circ \Delta _{-\rho }(F_{\cdot }(M_{\text{min}}))
\end{equation*}%
induces an isomorphism of homology groups, provided the left hand side has
finite-dimensional homology.

Therefore we are interested in calculating the homologies of 
\begin{equation*}
LT_{z}(\mathcal{I}(y,\mathbb{C}))
\end{equation*}%
for $z\in Y-Q_{\text{o}}$. Put $U=Q_{\text{o}}\cup Q_{\text{fp}}$. Thus $U$
is a Zariski open set isomorphic to $\mathbb{C}^{2}$. Let 
\begin{equation*}
j:Q_{\text{o}}^{\text{alg}}\rightarrow U^{\text{alg}}\text{ \ and \ }k:U^{%
\text{alg}}\rightarrow Y^{\text{alg}}
\end{equation*}%
denote the inclusions. Since $k_{\ast }(\mathcal{I}(y,\mathbb{C})\mid _{U^{%
\text{alg}}})\cong \mathcal{I}(y,\mathbb{C})$ and since $U$ is an affine
open set it follows from the base change that 
\begin{equation*}
LT_{z}(\mathcal{I}(y,\mathbb{C}))\cong 0\text{ \ for }z\in Y-U.
\end{equation*}%
On the other hand 
\begin{equation*}
\mathcal{I}(y,\mathbb{C})\mid _{U^{\text{alg}}}\cong j_{\ast }(\mathcal{O}%
_{Y^{\text{alg}}}\mid _{Q_{\text{o}}^{\text{alg}}})\text{.}
\end{equation*}%
Thus if $\left\{ z\right\} =Q_{\text{fp}}$ and $n>1$ then 
\begin{equation*}
LT_{z}(\mathcal{I}(y,\mathbb{C}))\cong \mathbb{C}
\end{equation*}%
since 
\begin{equation*}
j_{\ast }(\mathcal{O}_{Y^{\text{alg}}}\mid _{Q_{\text{o}}^{\text{alg}%
}})\cong \mathcal{O}_{Y^{\text{alg}}}\mid _{U^{\text{alg}}}
\end{equation*}%
Therefore 
\begin{equation*}
h_{p}(\Delta _{-\rho }(F_{\cdot }(M_{\text{min}})))\cong 
\begin{array}{c}
0\;\ \text{\ if }p\neq 0 \\ 
\left( \mathcal{O}_{Y}\mid _{S_{\text{o}}\cup S_{\text{fp}}}\right) ^{Y}%
\text{ \ if }p=0%
\end{array}%
\end{equation*}%
where $\left( \mathcal{O}_{Y}\mid _{S_{\text{o}}\cup S_{\text{fp}}}\right)
^{Y}$ denotes the extension by zero to $Y$ of the restriction of the sheaf
of holomorphic functions to $S_{\text{o}}\cup S_{\text{fp}}$. In particular, 
$\left( \mathcal{O}_{Y}\mid _{S_{\text{o}}\cup S_{\text{fp}}}\right) ^{Y}$
is the unique sheaf of DNF modules for the sheaf of holomorphic differential
operators on $Y$ whose sheaf cohomology vanishes in positive degrees and
whose global sections yield $M_{\text{min}}$.

Since $\mathcal{P}(y,\mathbb{C})\cong \mathcal{O}_{Y}\mid _{S_{\text{o}}}$,
we have the following short exact sequence 
\begin{equation*}
0\rightarrow \left( \mathcal{O}_{Y}\mid _{S_{\text{fp}}}\right)
^{Y}\rightarrow \left( \mathcal{O}_{Y}\mid _{S_{\text{o}}\cup S_{\text{fp}%
}}\right) ^{Y}\rightarrow \mathcal{P}(y,\mathbb{C})^{Y}\rightarrow 0.
\end{equation*}%
Taking global sections we obtain $M_{\text{min}}\cong \Gamma (S,$ $\mathcal{P%
}(y,\mathbb{C}))$.

\end{document}